\documentclass[12pt,leqno]{article}
 \def\thebibliography#1#2{\section*{#2
 }\list
 {[\arabic{enumi}]}{\settowidth\labelwidth{[#1]}\leftmargin\labelwidth
 \advance\leftmargin\labelsep
 \usecounter{enumi}}
 \def\newblock{\hskip .11em plus .33em minus .07em}
 \sloppy\clubpenalty4000\widowpenalty4000
 \sfcode`\.=1000\relax}
\newcommand{\binom}[2]{\left(\begin{array}{c}#1\\#2\end{array}\right)}
\renewcommand{\theequation}{\arabic{section}.\arabic{equation}}
\newcommand{\C}{{\bf C}}
\newcommand{\bC}{{\bf\bar C}}

\newcommand{\dpr}{{\partial}}
\newcommand{\varint}{{\frac{1}{2\pi}\int_{-\pi}^\pi}}
\newcommand{\U}{{\bf U}}
\newcommand{\e}{{e^{i\theta}}}
\title{Ahlfors' contribution to the theory of meromorphic functions}
\author{A. Eremenko\thanks{Supported by US-Israel Binational Science
Foundation and by NSF Grant DMS-9800084}}
\date{\today}

\begin{document}
\maketitle

\begin{abstract}
This is an expanded version of one of the
Lectures in memory of Lars Ahlfors in Haifa in 1996.
Some mistakes are corrected and references added. It contains a survey
of his work on meromorphic functions and related topics written
in 1929-1941.
\end{abstract}

This article is an exposition for non-specialists of Ahlfors' work
in the theory of meromorphic functions. When the domain is not specified
we mean meromorphic functions in the complex plane $\bf C$. 

The theory of meromorphic functions probably begins
with the book by Briot and Bouquet
\cite{BB} where the terms ``pole'', ``essential singularity'' and
``meromorphic'' were introduced and what is known now as
the Casorati--Weierstrass Theorem was stated for the first time.
A major discovery was Picard's theorem (1879) which says that a meromorphic
function omitting three values in the extended complex plane $\bC$ is
constant. The modern theory of meromorphic functions begins with attempts
to give an ``elementary proof'' of this theorem. These attempts culminated
in R. Nevanlinna's theory which was published first in 1925. Nevanlinna's
books \cite{N1} and \cite{N2}
were very influential and shaped 
much of
the research in function theory in 20th century. Nevanlinna theory, also
known as value distribution theory, was considered as one of the most important
areas of research in 1930-40, so it is not surprising that Ahlfors started
his career with work in this subject (besides he was Nevanlinna's student).
There are two very good sources about this early work of Ahlfors,
the first being Volume 1 of his collected papers \cite{CO},
where most of the papers have
been supplied with Ahlfors' own commentaries, the second Drasin's article
\cite{Dras}. I used both sources substantially in preparing this lecture,
trying to minimize intersections with Drasin's survey.
\vspace{.1in}
\setcounter{section}{1}
\setcounter{equation}{0}

\noindent
{\bf 1. Type problem}
\vspace{.1in}

All Ahlfors' papers on the theory of meromorphic functions are written in
the period 1929-1941 and they are unified by some general underlying ideas. 
The central problem is the problem of type. And the main method is the
Length -- Area
Principle in a broad sense. We start with explanation of
the type problem. 

According to the Uniformization theorem every open simply connected abstract
Riemann surface\footnote{Connected one dimensional complex analytic manifold}
$R$
is conformally equivalent to the complex plane $\C$ or to the unit disk $\U$.
In the first case $R$ is said to have {\em parabolic type} and in the second
case -- {\em hyperbolic type}. Now assume that a Riemann surface is given	
in some explicit way. How do we recognize its
type? This is the formulation of the type problem, which occupied a central
place in function theory in 1930-50. Ahlfors in
\cite[p. 84]{CO} gives credit to A. Speiser \cite{Speiser} for this formulation.

One very natural explicit way of describing a Riemann surface is this: let
$R$ be an abstract surface\footnote{Two real-dimensional
orientable triangulable topological
manifold} and let
\begin{equation}
\label{1}
\pi: R\rightarrow\bC
\end{equation}
be a topologically 
holomorphic map, which means that every point $x_0\in R$ has a
neighborhood $V$ and a local coordinate
$z:V\rightarrow\C$, $z(x_0)=0$ such that $f(x)=z^m(x)$ in $V$, where $m$ is a
positive integer. Points $x_0$ for which $m\neq 1$ are called critical points.
There is a unique conformal structure on $R$ which makes $\pi$ meromorphic.

The intuitive image connected with (\ref{1}) is the Riemann surface
``spread out over the sphere", which is an awkward English 
equivalent of ``\"{U}ber\-lagerungs\-fl\"{a}che''. 
We think of $R$ as a union of pieces of spheres (sheets) placed over
the sphere $\bC$ such that $\pi$ is the vertical projection. The sheets are
pasted together respecting the projection $\pi$.
This is how Riemann surfaces are usually introduced in elementary 
textbooks as Riemann surfaces
of multi-valued functions.

In what follows we reserve the term ``Riemann surface'' for a pair $(R,\pi)$
as in (\ref{1}). We shall frequently use {\em conformal metrics} on a Riemann
surface; these are Riemannian metrics compatible with the conformal structure.
In local coordinates such metric is expressed as $ds=p(z)|dz|$, where $p$
is a non-negative measurable function. In all situations we consider,
$p$ will be continuous and positive everywhere except at
 a set of isolated points.
A conformal metric can be always pulled back via a topologically holomorphic
map.
To a conformal metric corresponds an {\em area element} $d\rho=(i/2)p^2(z)
dz\wedge d{\bar z}=p^2(x+iy)dx\wedge dy$. The simplest example is the spherical Riemannian
metric. This is what one usually uses to measure the distance between
West Lafayette, IN and Haifa, for example.

The type problem in this setting is to find out from the geometric information
about $(\ref{1})$ whether $R$ is conformally equivalent to $\C$ or to $\U$
(assuming that $R$ is open and simply connected).
For example Picard's theorem gives the following necessary condition of
parabolic type: if $R$ is parabolic then it covers all points of $\bC$ with
at most two exceptions.

In \cite{A0} Ahlfors gives a sufficient condition for parabolic type.
Assume that $\pi$ has the following property: every curve in $\bC$ has
a lifting to $R$.
This means that the only singularities of the multiple valued analytic
function $\pi^{-1}$ are algebraic branch points (critical points of $\pi$).
 We choose a point $x_0\in R$
and exhaust $R$ with disks (with respect to the pullback of
the spherical metric) of radius $r$ centered at $x_0$. Let $n(r)$ be the number
of critical points of $\pi$ in the disk of radius $r$. If
\begin{equation}
\label{2}
\int^\infty\frac{r\,dr}{n(r)}=\infty
\end{equation}
then $R$ is of parabolic type. If for example $R=\C$ and $\pi$ is an elliptic
function then the pullback of the spherical metric differs from the
Euclidean metric by a factor bounded from above and below, so
$n(r)/r^2$ is bounded from above and below and we see that condition
(\ref{2}) gives the right order of growth for $n(r)$.

In \cite{A1} Ahlfors
goes further and
considers {\em arbitrary} conformal Riemannian metrics on $R$ such that
$R$ is complete. Let $L(r)$
be the length of the circumference of the circle of radius $r$ centered at $z_0$
with respect to such a metric. Then $R$ has parabolic type {\em if
and only if there
exists a metric of the described type with the property}
\begin{equation}
\label{la}
\int^\infty\frac{dr}{L(r)}=\infty. 
\end{equation}
Here is the proof of (\ref{la}) which demonstrates how the Length --
Area principle works. 
We may identify $R$ with a disk in $\C$,
$R=\{ z\in\C:|z|<r_0\},\; r_0\leq\infty$, 
and write the length element of a conformal
metric as $ds=p(z)|dz|$.
Let $\Gamma_t$ be the circumference of the
disk of radius $t$, centered at $0$ with respect
to $ds$.
Then
$$2\pi\leq\int_{\Gamma_t}|d\log z|=\int_{\Gamma_t}
\frac{|d\log z|}{ds}ds,$$ and by the Schwarz inequality
$$4\pi^2\leq\int_{\Gamma_t}ds\int_{\Gamma_t}
\left(\frac{|d\log z|}{ds}\right)^2ds=L(t)\int_{\Gamma_t}
\frac{|d\log z|^2}{ds^2}ds.$$
So
$$4\pi^2\int_{r'}^{r''}\frac{dt}{L(t)}\leq
\int_{r'}^{r''}dt\int_{\Gamma_t}\left(\frac{|d\log z|}{ds}\right)^2 ds.$$
The last double integral is the area swept by $\Gamma_t$
as $t$ varies from $r'$ to $r''$,
measured with respect to the conformal metric $|d\log z|=
|z|^{-1}|dz|$.
If our Riemann surface is hyperbolic,
that is $r_0<\infty$, this area stays
bounded as $r''\to r_0$ and the integral
(\ref{la}) converges for every conformal
metric
$ds$. If the surface is parabolic we just take $ds=|dz|$.

Of course (\ref{la}) does not give any effective solution of the type problem,
but almost all known criteria of parabolic type, including (\ref{2}), where
the spherical metric is used, can be obtained by a special
choice of concrete metric on $R$.  
Two sources about later developments in the type problem are
\cite{Volkovyski} and \cite{Wittich}. 

Assume now that the type of $R$ is known, and let $\theta$ be a conformal
map of the disk or plane $D=\{ z:|z|<r_0\},\; r_0\leq\infty$ onto $R$.
Then $f=\pi\circ\theta$ becomes
a meromorphic function in $D$ which is almost uniquely defined (up to
a conformal automorphism of $D$). One may be interested in how the geometric
properties of $\pi$ are connected with the asymptotic behavior of $f$.
By asymptotic behavior of $f$ we mean for example how the number of solutions
of
equations $f(z)=a$ in the disk $\{ z: |z|\leq r\}$ behaves when $r\rightarrow r_0$. This is the subject of the value distribution theory.
\newpage 

\noindent
{\bf 2. Value distribution theory}
\setcounter{section}{2}
\setcounter{equation}{0}
\vspace{.1in}

If $f$ is a rational function its value distribution is controlled by
its degree $d$, which is the number of preimages of a generic point.
The main tool
of the value distribution theory of meromorphic functions of R. Nevanlinna
is the characteristic function $T_f(r)$ which replaces the degree in the case
when $f$ is transcendental. We explain the version of the Nevanlinna theory
which was found by Shimizu and Ahlfors independently of each other.

Let us denote by $n(r,a)=n_f(r,a)$ the number of
solutions of the equation $f(z)=a$
in the disk $\{ z: |z|\leq r\}$, counting multiplicity. Here $a\in \bC$.
By the Argument Principle and the Cauchy--Riemann equations we have
$$n(r,a)-n(r,\infty)=\frac{1}{2\pi i}\int_{|z|=r}\frac{f'}{f-a}dz=
\frac{r}{2\pi}\frac{d}{dr}\int_{-\pi}^\pi\log |f(r\e)-a|\,d\theta.$$
We divide by $r$ and integrate with respect to $r$, assuming for a moment
that $f(0)\neq a,\infty$ and using the notation\footnotemark
\footnotetext{If $f(0)=a$
this has to be regularized in the following way:
$$N(r,a)=\int_0^r\{ n(t,a)-n(0,a)\} t^{-1}dt+n(0,a)\log r.$$}
$$N(r,a)=\int_0^r\frac{n(t,a)}{t}dt\quad$$
to obtain
\begin{equation}
\label{3}
\varint\log |f(r\e)-a|\, d\theta=\log |f(0)-a|+N(r,a)-N(r,\infty).
\end{equation}
This is the Jensen formula.

Let $d\rho=p^2(u+iv)du\wedge dv$ be the area element of a conformal metric, 
normalized such that the area of $\bC$ equals $1$. We
integrate (\ref{3}) with respect to $a\in \bC$ against $d\rho$
and obtain
\begin{equation}
\label{3.2}
\varint U(f(r\e))\, d\theta=U(f(0))+\int_\bC N(r,a)\, d\rho(a)-N(r,\infty),
\end{equation}
where
\begin{equation}
\label{3.5}
U(w)=\int_{\bC}\log |w-a|\,d\rho(a).
\end{equation}
If we choose $d\rho$ to be the normalized spherical area element that is
$$p(w)=\frac{1}{\sqrt{\pi}(1+|w|^2)}$$ then changing the order of integration
shows
\begin{equation}
\label{4}
T(r):=\int_\bC N(r,a)\,d\rho(a)=\int_0^r\frac{A(t)\, dt}{t},
\end{equation}
where
$$A(r)=\frac{1}{\pi}\int_{|z|\leq r}
\frac{|f'|^2}{(1+|f|^2)^2}dm,\quad\mbox{where}\; dm\;\mbox{is the Euclidean 
area element in}\;\C.$$
The geometric interpretation of $A(r)$ is the area of the disk
$|z|\leq r$ with respect to the pullback of the spherical metric,
or in other words, the
{\em average covering number} of the Riemann sphere by the restriction of
$f$ to the disk $|z|\leq r$. The function $T(r)$ defined in (\ref{4}) is called
the
{\em characteristic} of $f$. It is analogous to the degree of a rational function
in the sense that it measures the number of preimages of a generic
point.\footnote{In fact $T(r)$ also has algebraic properties of degree:
it is a {\em logarithmic height} in the field of meromorphic function.
The importance of this fact was fully recognized only recently
\cite{Vojta,Lang,Eremenko}} 
If $f$ is rational we have $T_f(r)=\deg f\cdot\log r+O(1)$, and for transcendental
$f$ we always have $T(r)/\log r\rightarrow\infty$ as $r\to\infty$.

Now if we use the spherical area element $d\rho$ in (\ref{3.5}), the integrals can be evaluated
which gives $U(w)=\log\sqrt{1+|w|^2}=\log ([w,\infty])^{-1}$, where $[,]$
stands for the {\em chordal} distance on the Riemann sphere.
Thus the first term in (\ref{3.2}) is
$$m(r,\infty):=\varint\log\frac{1}{[f(r\e),\infty]}d\theta.$$
and in general we can define
$$m(r,a):=\varint\log\frac{1}{[f(r\e),a]}d\theta,\quad a\in\bC.$$
This is called the {\em proximity function}\footnote{Nevanlinna's original definition of proximity function differs by an additive term which is bounded when 
$r\to\infty$.}. It becomes large when $f$ is close
to $a$ in the average on the circle $|z|=r$. It is important that the
proximity function is non-negative.  With these notations (\ref{3.2}) can
be rewritten as
$$N(r,\infty)+m(r,\infty)=T(r)+m(0,\infty).$$
Now we notice that $T(r)$ does not change if we replace $f$ by $L\circ f$ where
$L$ is a rotation of the sphere, which is a conformal automorphism which
preserves spherical distance. Thus we obtain
\begin{equation}
\label{FMT}
\; N(r,a)+m(r,a)=T(r)+m(0,a)=T(r)+O(1), \;r\rightarrow\infty,\; a\in\bC.
\end{equation}
This is the First Main Theorem (FMT) of Nevanlinna
in the form of Shimizu--Ahlfors. It implies by the way that
\begin{equation}
\label{FMT-bis}
N(r,a)\leq T(r)+O(1),\;r\to\infty,\;\mbox{for every}\; a\in\bC,
\end{equation}
because the proximity function is
non-negative. When compared with (\ref{4}) this shows that the points
$a$ for which $N(r,a)$ is substantially less then $T(r)$ should be exceptional.
This was put in a very precise form by Ahlfors and Nevanlinna \cite{N2},
who extended
the earlier results by Valiron and Littlewood.
We state
only one theorem of this sort:
\vspace{.1in}

\noindent\addtocounter{equation}{+1}
{(\theequation)\em For every $\epsilon>0$ we have $m(r,a)=O(T(r))^{1/2+\epsilon}$ for all
$a\in\bC\backslash E$, where $E$ is a set of zero logarithmic capacity.}
\vspace{.1in}

\noindent
The description of the exceptional set in this theorem cannot be
substantially improved without making any additional assumptions about $f$.
Very subtle
examples constructed by Hayman \cite{Hayman} show that
\vspace{.1in}

\noindent\addtocounter{equation}{+1}
(\theequation){\em Given an
$F_\sigma$ set $E\subset\bC$ of zero logarithmic capacity there exists a meromorphic
function $f$ and a sequence $r_k\to\infty$ such that $m_f(r_k,a)\sim T_f(r_k)$
for all $a\in E$.}
\vspace{.1in}

\noindent
On the other hand, a lot of work has been done to improve
the exceptional set for various subclasses of meromorphic functions
(see \cite{Dras} for a survey of these results).

The functions $T_f(r),N_f(r,a)$ and $m_f(r,a)$ provide a convenient way to
describe the asymptotic behavior of $f$.
The characteristic $T(r)$ measures the
growth of the number of solutions of $f(z)=a$ in the disks $|z|\leq r$.
The First Main Theorem shows that $T(r)$ provides upper bound for
this number {\em for all $a\in\bC$} and (\ref{4}) or (2.7) 
show that for typical $a$ this upper bound is attained.
It is natural to define the {\em order} of a meromorphic function $f$ by
$$\lambda=\limsup_{r\to\infty}\frac{\log T(r)}{\log r}.$$
\vspace{.1in}

\noindent
{\bf 3. Simply connected parabolic Riemann surfaces}.
\setcounter{section}{3}
\setcounter{equation}{0}
\vspace{.1in}

There are few known sufficient conditions for parabolic type. The one
mentioned in Section 1 requires that all paths can be lifted via $\pi$
so it cannot be applied to Riemann surfaces (\ref{1}) with ``singularities''.
Let us define this notion precisely. 

Let $D(a,r)\subset\bC$ be the disk with respect to the
spherical metric of radius $r$ centered at $a\in\bC$. We fix $a\in\bC$ and
consider a function $\sigma$ which to every $r>0$ puts into correspondence
a component $\sigma(r)$ of the preimage $\pi^{-1}(D(a,r))$ in such a way that
$r_1>r_2$ implies $\sigma(r_1)\subset\sigma(r_2)$. Now there are two
possibilities:
\newline
a) $\cap_{r>0}\sigma(r)=\{\mbox{one point}\}\subset R$ or
\newline
b) $\cap_{r>0}\sigma(r)=\emptyset$.
\newline
In the case b) we say that $\sigma$ defines a (transcendental)
{\em singular point} over $a$
of the
Riemann surface $(R,\pi)$. Thus all possible $\sigma$'s are in one to one
correspondence with all ``points'', ordinary and singular.

A singular point $\sigma$ over $a$ is called {\em logarithmic} if for some
$r>0$ the restriction $\pi|_{\sigma(r)}$ is a universal covering over
$D(a,r)\backslash\{ a\}$.

The following classification belongs to Iversen: a singular point $\sigma$
over $a$ is called {\em direct} if there exists $r>0$ such that $\pi$ omits the
value $a$ in $\sigma(r)$. If such $r$ does not exist, $\sigma$ is called
indirect. So logarithmic singular points are direct.

\noindent
{\em Examples}. Let $R=\C$. If $\pi=\exp$ there are two singular points:
one over $0$ and one over $\infty$. Both are direct (logarithmic).
If $\pi(z)=\sin z/z$
there are two indirect singular points over $0$ and two direct (logarithmic)
singular points
over infinity. If $\phi(z)=z\sin z$ then there is one direct singular point over
infinity, which is not logarithmic.

Nevanlinna \cite{Nev} proved that an open simply connected Riemann surface having
no critical points and finitely many singular points is of parabolic
type. This was extended by Elfving \cite{Elfving} to allow finitely many
critical points. In fact they obtained much more precise information about meromorphic functions
associated  with these surfaces:
\vspace{.1in}

\addtocounter{equation}{+1}
\noindent
{(\theequation) \bf Theorem}. {\em Let $\pi:R\to\bC$ be a simply connected
Riemann
surface,
$\theta$ a conformal mapping of the disk $|z|<r_0\leq\infty$ onto $R$,
and $f=\pi\circ\theta$. Assume that $\pi$ has only finitely many
critical points
and that $(R,\pi)$ has a finite number $p\geq 2$
of singular
points. Let $d(a)$ be the number of singular
points over $a$. Then:
\newline
$(i)$ $(R,\pi)$ is of parabolic type, so that $r_0=\infty$;
\newline
$(ii)$ $T_f(r)\sim cr^{p/2},\; r\to\infty$, where $c>0$ is a constant, and
\newline
$(iii)$ $m_f(r,a)\sim (2cd(a)/p)r^{p/2},\; r\to\infty,\; a\in\bC$.}
\vspace{.1in}

The simplest example of the situation described by this theorem is
$p=2$ and no critical points, then $f$ is an exponential. If $p=3$ and there are no critical
points, $f$ can be expressed in terms
of the Airy functions which satisfy the differential equation $y''+zy=0$.
Nevanlinna's proof is based on asymptotic
integration of certain linear differential equations with polynomial
coefficients, which makes it impossible
to extend the method to any substantially wider class of Riemann surfaces
(at least nobody has ever succeeded in this).
In the same journal issue where Nevanlinna's theorem (3.1) appeared,
Ahlfors' paper \cite{A3} was published, which gave a completely different
approach to the problem. It is Ahlfors' approach which was the base of all
subsequent generalizations. Roughly speaking, the argument goes in the following
way. One dissects $R$ in (\ref{1}) into pieces such that for each piece
an explicit conformal map of a plane domain onto this piece
can be found. Then one pastes
these pieces together and obtains a homeomorphism of $\C$ onto $R$.
The point is to perform this construction in such a way that the explicitly
constructed
homeomorphism $\theta_1:\C\to R$ 
is as close to conformal as possible. If one finds a quasiconformal
homeomorphism $\theta_1$ this implies that $R$ has parabolic type.
One can go further and find $\theta_1$ 
close to a conformal $\theta$ in the following sense:
$$\theta=\theta_1\circ\phi,\quad\mbox{where}\quad\phi(z)\sim cz,\; z\to\infty.$$
Then one can derive the approximate formulas for $N(r,a)$ if $\theta_1$ is
described explicitly. (Ahlfors' argument was in fact more complicated, but this
is what was eventually distilled from his paper). 
The distortion theorem used 
to show that $\theta$ is close to $\theta_1$ is called the
Teichm\"{u}ller--Wittich--Belinskii Theorem: {\em if $\phi$ is a quasiconformal
homeomorphism of the plane and
$$\int_{|z|\geq 1} k(z)\frac{dm}{|z|^2}<\infty,\quad\mbox{where}\quad
k(z)=|\phi_{\bar z}/\phi_z|$$
then $\phi$ is conformal at infinity, that is the limit 
$\lim_{z\to\infty}\phi(z)/z\neq 0$ exists.
}

Using the described approach more and more general classes of
open simply connected
Riemann surfaces were subsequently introduced, their type determined,
 and the value distribution of
corresponding meromorphic functions studied. We mention the work of
K\"{u}nzi, Wittich, Schubart, P\"{o}schl, Ullrich, Le Van Thiem, Huckemann, 
H\"{a}llstr\"{o}m, and Goldberg. Most of the results are
described 
in the books of Wittich \cite{Wittich} and Goldberg-Ostrovskii \cite{GO} and
in Goldberg's paper \cite{piends}. For further development 
see \cite{DW} -- the only paper on the subject written in English.

Let us assume now that $(R,\pi)$ is a Riemann surface of parabolic type and see how
the topological properties of $\pi$ 
influence the asymptotic behavior of the meromorphic
function $f$. 
The most famous result in this direction is the
Denjoy--Carleman--Ahlfors Theorem \cite{A2}. 
\vspace{.1in}

\noindent\addtocounter{equation}{+1}
{(\theequation) \bf Theorem.} {\em If a simply connected parabolic Riemann surface has $p\geq 2$
direct singular points, then the corresponding meromorphic function satisfies
$$\liminf_{r\rightarrow\infty}r^{-p/2}T(r)>0.$$
Thus the number of direct singularities is at most $\max\{1,2\lambda\}$,
where $\lambda$ is the order of the meromorphic function.} 
\vspace{.1in}

\noindent
The story of this theorem is explained in detail in M. Sodin's talk
\cite{Sodin}. We add only few remarks. Heins \cite{Heins} proved that
\vspace{.1in}

\noindent\addtocounter{equation}{+1}
(\theequation) {\em
the set of projections of direct singular points of a parabolic Riemann surface
is at most countable.}
\vspace{.1in}

\noindent
On the other hand, simple examples show that the set of direct singular points
itself may have the power of the continuum.
One cannot say anything about the size of the set of projections of
{\em all} singular
points, even if the growth of $T_f$ is restricted: for every $\lambda\geq 0$ there exist meromorphic functions
of order $\lambda$ such that every point in $\bC$ has a singular point
over it 
\cite{Erem}. 
Goldberg \cite{Goldberg}
generalized Ahlfors' theorem by including a certain
subclass of indirect singular points, which are called $K$-points.
The property of being a $K$-point depends only on the restriction of $\pi$
to $\sigma(r)$ with arbitrarily small $r>0$, but only for some narrow
classes of Riemann surfaces is it known how to determine effectively whether
a singular point is a $K$-point.
\vspace{.1in}

\noindent
{\bf 4. The Second Main Theorem of the value distribution theory}
\setcounter{section}{4}
\setcounter{equation}{0}
\vspace{.1in}

To formulate the main result of the value distribution theory we 
denote by $n_1(r)=n_{1,f}(r)$
the number of critical
points of meromorphic function $f$ in the disk $|z|\leq r$, counting
multiplicity. It is easy to check that
\begin{equation}
\label{ram}
n_{1,f}(r)=n_{f'}(r,0)+2n_f(r,\infty)-n_{f'}(r,\infty).
\end{equation}
Now we apply the averaging as above:
$$N_1(r)=N_{1,f}(r):=\int_0^r\frac{n_1(t)\,dt}{t}.$$
If $0$ is a critical point, the same regularization as before has to be made. 
The Second Main Theorem (SMT) says that {\em for every finite set $\{ a_1,\ldots,a_q\}
\subset\bC$ we have
\begin{equation}
\label{SMT}
\sum_{j=1}^q m(r,a_j)+N_1(r)\leq 2 T(r)+S(r),
\end{equation}
where $S(r)=S_f(r)$ is a small ``error term'', $S_f(r)=o(T(r))$
when $r\to\infty, r\notin E$, where $E\subset[0,\infty)$
is a set of finite measure}\footnote{In fact $S(r)$ has much better
estimate. Recently there was a substantial
activity in the study of the best possible estimate of this error term,
see for example \cite{Hink}. On the other hand Hayman's 
examples (2.8) 
show that in general the error term may not be $o(T(r))$
for all $r$, so an exceptional set $E$ is really required.}.
The SMT may be regarded as a very precise way of saying that the term
$m(r,a)$ in the FMT (\ref{FMT}) is relatively small for most $a\in \bC$.
It is instructive to rewrite (\ref{SMT}) using (\ref{FMT}) in the
following form. Let ${\bar N}(r,a)$ be the averaged counting function
of {\em different} solutions of $f(z)=a$, that is this time we don't
count multiplicity. Then $\sum_jN(r,a_j)\leq\sum_j {\bar N}(r,a_j)+N_1(r)$
and we obtain
\begin{equation}
\label{SMT-bis}
\sum_{j=1}^q{\bar N}(r,a_j)\geq (q-2)T(r)+S(r).
\end{equation}
Now Picard's theorem is an immediate consequence: if three values
$a_1,a_2$ and $a_3$ are omitted by a meromorphic function $f$, then 
$N_f(r,a_j)\equiv 0,\; 1\leq j\leq 3$, so the left side of (\ref{SMT-bis})
is zero
and we obtain $T_f(r)=S_f(r)$, which implies that $f$ is constant.
Here is a more refined
\vspace{.1in}

\addtocounter{equation}{+1}
\noindent
{(\theequation)\bf Corollary from the SMT.} {\em Let $a_1,\ldots, a_5$ be
five points on the Riemann sphere. Then at least one of the equations
$f(z)=a_j$ has simple solutions.}
\vspace{.1in}

Indeed, if all five equations
have only multiple solutions then $N_1(r,f)\geq (1/2)\sum_{j=1}^5 N(r,a_j).$
When we combine this inequality with SMT (\ref{SMT}) it implies
$(5/2)T(r)\leq 2T(r)+S(r)$, so $f=\mbox{const}$. 

For most ``reasonable'' functions, like Nevanlinna's
functions described in Theorem (3.1),
the SMT tends to be an asymptotic equality
rather then inequality; the most general class of functions for which
this is known consists of meromorphic functions whose
critical and singular points lie over a finite set. This is due to
Teichm\"{u}ller; 
see, for example \cite[Ch. 4]{Wittich}).

The purpose of Ahlfors' paper \cite{A6},
as he explains in the commentary, ``was to derive the main results of
Nevanlinna theory of value distribution in the simplest way I knew how.''
The proof presented in \cite{A6} was the source of most generalizations
of Nevanlinna theory to higher dimensions. We will return to these
generalizations in Section 8 and now Ahlfors' proof will be presented.
A reader not interested in proofs may skip to the end of this section.
\vspace{.1in}

\noindent
{\em Ahlfors' proof of SMT}. 
Let $d\omega$ be the spherical area element, so that
$$d\omega=\frac{dudv}{\pi (1+|w|^2)^2},\quad w=u+iv.$$
 
We consider another area element on the
Riemann sphere, $d\rho=p^2d\omega$, where $p$ is 
given by
\begin{equation}
\label{sing}
\log p(w):=\sum_{j=1}^q\log\frac{1}{[w,a_j]}-2\log\left(\sum_{j=1}^q\log
\frac{1}{[w,a_j]}\right)+C,
\end{equation}
where $[,]$ is the chordal distance,
and $C>0$ is chosen so that 
$$\int_\bC d\rho=1.$$
(The sole purpose of the second term in the definition
of $p$ in (\ref{sing}) is to make this integral
converge, without altering much the behavior near $a_j$ which is determined by
the first term). 

We pull back this $d\rho$ via $f$ and write the change of the variable formula:
\begin{equation}
\label{fund}
\int_\bC n(r,a)d\rho(a)=\int_0^r\int_{-\pi}^\pi p^2(w)\frac{|w'|^2}{(1+|w|^2)^2}
t\,d\theta dt,\quad w=f(te^{i\theta}).
\end{equation}
Now we consider the derivative of the last double integral with respect to $r$,
divided by $2\pi r$:
$$\lambda(r):=\varint\frac{|w'|^2}{(1+|w|^2)^2}p^2(w)\,d\theta,\quad w=f(re^{i\theta}).$$
Using the integral form of the arithmetic-geometric means 
inequality\footnote{
$\quad\frac{1}{b-a}\int_a^b\log g(x)\,dx\leq\log\left\{\frac{1}{b-a}\int_a^bg(x)\,dx
\right\}$
} we obtain
\begin{equation}
\label{prep}
\;\;\log\lambda(r)\geq\frac{1}{\pi}\int\log p(w)\, d\theta-
\frac{1}{\pi}\int \log (1+|w|^2)d\theta+\frac{1}{\pi }\int\log |w'|\, d\theta.
\end{equation}
The first integral in the right-hand side of (\ref{prep})
is approximately evaluated using (\ref{sing})
(the second summand in (\ref{sing})
becomes irrelevant because of another $\log$);
the second integral equals
$4m(r,\infty)$ and the third one
 is evaluated using Jensen's formula (\ref{3}):
$$(2+o(1))\sum_{j=1}^qm(r,a_j)+
2\left\{ N(r,0,f')-N(r,\infty,f')-2m(r,\infty)\right\}
\leq\log\lambda(r)$$
The expression inside the brackets is equal to $N_1(r)-2T(r)$ (by 
definition of $N_1$ and the FMT (\ref{FMT}) applied with $a=\infty$), so 
\begin{equation}
\label{prelim}
\sum_{j=1}^qm(r,a_j)+N_1(r)-(2+o(1))T(r)\leq\frac{1}{2}\log\lambda(r).
\end{equation}
To estimate $\lambda$ we return to the left side of (\ref{fund}).
Integrating twice and using the FMT we obtain
$\int_0^r\frac{dt}{t}\int_0^t\lambda(s)sds=\int_\bC N(r,a)\,d\rho(a)\leq
T(r)+O(1).$
Now the argument is concluded with the following elementary calculus lemma:
\vspace{.1in}

\noindent
{\em If $g$ is an increasing function on $[0,\infty)$,
tending to infinity, then 
$g'(x)\leq g^{1+\epsilon}(x)$ for all $x\notin E$, where $E$ is a set of
finite measure.}
\vspace{.1in}

\noindent
Applying this lemma twice we conclude that
$\log\lambda(r)=S(r)$ which proves the theorem.
\vspace{.2in}

\noindent
{\bf 5. Ahlfors'  \"{U}berlagerungsfl\"{a}chentheorie}
\setcounter{section}{5}
\setcounter{equation}{0}
\vspace{.1in}

In a series of papers written in 1932-33 Ahlfors developed
a different approach to
value distribution theory, which is based on a generalization of the
Riemann--Hurwitz formula.

To explain his motivation we return to the type problem. As was explained
in Section 1, Picard's theorem may be considered as a necessary condition of
parabolic type. But this condition is very unstable: a small perturbation
of $(R,\pi)$ destroys the property that a value is omitted. The same
can be said about (4.4) which may also be regarded as a necessary condition for 
parabolic type. This instability was noticed by A. Bloch \cite{Bloch},
who stated the 
following general Continuity Principle, which I prefer to cite from
Ahlfors' survey 
\cite{A4} {\em une proposition de nature qualitative, exacte avec
un certain \'{e}nonc\'{e}, demeure encore exacte si l'on modifie
les donn\'{e}es de la proposition, en leur faisant subir une
d\'{e}formation continue}. In accordance with this principle Bloch
conjectured among other things the following:
\vspace{.1in}

\noindent\addtocounter{equation}{+1}
(\theequation)
 {\em Let $D_1,\ldots,D_5$ be five Jordan
domains with disjoint closures, and $f$ be a non-constant
meromorphic function
in $\C$. Then there is a bounded simply connected domain $D\subset \C$
which is mapped by $f$ homeomorphically onto one of the domains $D_j$.}
\vspace{.1in}

\noindent
This conjecture improves (4.4) in accordance with the Continuity Principle.
A similar improvement of Picard's theorem would be that 
\vspace{.1in}

\noindent
{\em Given three
disjoint domains and a non-constant meromorphic function $f$,
there is at least
one bounded component of the $f$-preimage  of one of these domains.}
\vspace{.1in}

In \cite{A6} Ahlfors proved all these conjectures. His technically hard proof
was a combination of topological considerations with sophisticated
distortion theorems. Finally he developed a general theory \cite{A7} which probably
constitutes his most original contribution to the study of meromorphic
functions. The theory has a metric-topological nature; all
complex analysis in it is reduced to one simple application of the Length -- Area
Principle. This makes the theory flexible enough to treat quasiconformal
mappings because for such mappings the Length -- Area
argument also works. In fact the word ``quasiconformal'' was for the first time
used in this paper. In his commentary in \cite{CO} to \cite{A7} 
Ahlfors writes: 
\newline
``I included this more general situation in my paper but with pangs of
conscience because I considered it rather cheap padding... Little did I know at
that time what an important role quasiconformal mappings would come to play in my own work''.

We start with several definitions and elementary facts from the topology of
surfaces.
A bordered surface of finite type
is a closed region on a compact orientable surface bounded by finitely many
simple closed curves.
(Compact orientable surfaces are spheres with finitely
many handles attached). The Euler characteristic of a bordered surface
$\chi$ is\footnote{We follow Ahlfors' notations. In modern literature
the Euler characteristic is usually defined with the opposite sign.} 
$2g-2+k$, where $g$ is the number of handles (genus) and $k$ is the number
of holes (boundary curves).
We consider a topologically holomorphic map
of two bordered surfaces of finite type, $\pi:R\to R_0$. If $\pi(\dpr R)\subset
\dpr R_0$ then $\pi$ is a ramified covering. In this case we have
\newline
(i) there is a number $d$, called the degree of $\pi$ such that every point in $R_0$
has $d$ preimages, counting multiplicity;
\newline
(ii) $\chi(R)\geq d\chi(R_0)$.
\newline
In fact a more precise relation holds, $\chi(R)=d\chi(R_0)+i(\pi)$, where
$i(\pi)$ is the number of critical points of $\pi$, counting multiplicity.
This is the Riemann--Hurwitz formula.

Ahlfors' theory extends these two facts to the case when $\pi$ is {\em almost}
a ramified covering, that is the part of the boundary $\dpr_0 R\subset\dpr R$ which is mapped to the
interior of $R_0$ is relatively small in the sense we are going to define now.
This part $\dpr_0 R$ is called the relative boundary.

Let $R_0$ be a region in a compact surface equipped 
with a Riemannian metric. A curve is called regular if it is piecewise smooth.
A region is called regular if it is bounded by finitely many piecewise smooth
curves.\footnote{Ahlfors defines the precise degree of regularity for which his
theory works. Here we make stronger smoothness assumptions only for simplicity.}
We assume that $R_0$ is a regular region. The restriction of the Riemannian metric 
to $R_0$ is called $\rho_0$ and its pullback to $R$ is called $\rho$.
 
Let $L$ be the length of the relative boundary $\dpr_0 R$. Everything on $R$
is measured using the pulled back metric $\rho$. We use the symbol $|.|$ for
the area of a regular region or for the length of a regular curve, depending
on context. We define the {\em average number of sheets} of $\pi$ by
$S:=|R|/|R_0|$. Similarly, if $X_0\subset R_0$ is a regular region or curve,
we define the average number of sheets over $X_0$ as $S(X_0)=|X|/|X_0|$,
where $X=\pi^{-1}(X_0).$ 

The following results are the First and Second Main Theorems
of Ahlfors' theory.
\vspace{.1in}

\noindent\addtocounter{equation}{+1}
(\theequation)
{\em For every regular region or curve $X_0\subset R_0$ there exists a constant
$k$ depending only of $R_0$ and $X_0$ such that}
$$|S(X_0)-S|<kL.$$
\addtocounter{equation}{+1}
(\theequation)
{\em There exists a constant $k$ depending only of $R_0$ such that}
$$\max\{\chi(R), 0\}\geq S\chi(R_0)-kL.$$

It is important in these theorems that $k$ does not depend on $R$ and $\pi$.
They contain useful information if $L$ is much smaller then $S$, when they
imply that
properties (i) and (ii) above hold approximately.

The proofs of (5.2) and especially (5.3)
are rather sophisticated\footnote{Ahlfors won one
of the first two Fields Medals for this in 1936}, though elementary. The proof
of (5.3) was later simplified by Y. Toki \cite{Toki}.
There is another recent simplified proof by de Thelin
\cite{Th}.

Now we explain how this applies to open Riemann surfaces
and meromorphic functions. Let us consider a meromorphic function
$$f:R\rightarrow \bC,$$
where $R$ is the plane $\C$ equipped with the pullback
of the spherical metric $\rho_0$ on $\bC$.
One can consider any compact Riemann surface  with any Riemannian metric instead
of $\bC$. A surface with a Riemannian metric $R$ is called {\em regularly exhaustible} if
there exists an increasing sequence of compact regular subregions
$F_1\subset F_2\subset\ldots,\;\cup F_j=R$ with the property $|\dpr F_j|/|F_j|\to 0$
as $j\rightarrow\infty$. It is easy to prove using the Length -- Area
Principle that $\C$ is always regularly exhaustible, no matter what the
conformal metric is, and a sequence of concentric Euclidean disks can be always taken as
$F_j$. Thus one can apply (5.2) and (5.3) to the restrictions of a meromorphic
function $f$ on the disks $|z|\leq r$. Then the average number of sheets
$S=A_f(r)$ is the same $A(r)$ as in (\ref{4}).
We denote the length of the relative boundary by $L(r)$, this is nothing but
spherical length of the $f$-image of the circle $|z|=r$.
The Euler characteristic of the
sphere is negative, $\chi(\bC)=-2$, so (5.3) gives nothing. To extract useful
information from (5.3) we consider a region $R_0\subset\bC$ obtained by
removing from $\bC$ a collection of $q\geq 3$ disjoint spherical
disks $D_j$, so that $\chi(R_0)=q-2$. Now we apply (5.3) to the restriction
of $f$ on $R(r):=f^{-1}(R_0)\cap\{ z:|z|\leq r\}$. The Euler characteristic of
this region is the number of boundary curves minus two. A component
of $f^{-1}(D_j)\cap\{ z:|z|\leq r\}$
is called an {\em island} if it is relatively compact in
$|z|<r$, and  all other components are called {\em peninsulae}.
If ${\bar n}(r,D_j)$ is the number of islands over $D_j$, then 
$\chi(R(r))=\sum_j{\bar n}(r,D_j)-1$. Thus from (5.2) and (5.3) we obtain
\begin{equation}
\label{sch}
\sum_{j=1}^q {\bar n}(r,D_j)\geq (q-2)A(r)+o(A(r))
\end{equation}
for those $r$ for which $L(r)/A(r)\to 0$.
This result is called the Scheibensatz. It has to be compared with the
SMT of Nevanlinna in the form (\ref{SMT-bis}). First of all, (\ref{sch})
can be integrated with respect to $r$ to obtain
\begin{equation}
\label{N}
\sum_{j=1}^q {\bar N}(r,D_j)\geq (q-2)T(r)+o(T(r))
\end{equation}
with a somewhat worse exceptional set for $r$ than in (\ref{SMT-bis})
\cite{Miles}, \cite{ES}.
Now if $a_j$ lie in the interiors of $D_j$ for $j=1,\ldots,q$ then
evidently ${\bar N}(r,a_j)\geq {\bar N}(r,D_j)$ because the restriction of
$f$ on an island is a ramified
covering over $D_j$. Thus (\ref{N}) is better than (\ref{SMT-bis})
because (5.5) does not count $a$-points in peninsulae.

The Five Islands Theorem (5.1) follows from (\ref{sch}) in the same way
as (4.4) follows from (4.3).

A substantial complement to this theory was made by Dufresnoy
\cite{Duf}, who invented the way of deriving normality criteria
from Ahlfors' theory using a simple argument based on an isoperimetric
inequality (see also \cite{Haymanbook} for a very clear and self-contained
exposition of the whole theory and its applications).
The normality criteria we are talking about are expressions of another
heuristic principle\footnote{several rigorous results which may be considered
as implementations of this principle are known: \cite{Zalc,Lang}} 
sometimes associated with the name of Bloch
(based on his paper \cite{Bloch}) 
\vspace{.1in}

\noindent\addtocounter{equation}{+1}
(\theequation)
{\em that every condition which implies
that a function meromorphic in $\C$ is constant, when applied to a family
of functions in a disk, should imply that this family is normal, preferably
with explicit estimates}.
\vspace{.1in}

For example:
\vspace{.1in}

\noindent
{\em Given five regular Jordan regions with disjoint closures, consider
the class of meromorphic functions in a disk which have no simple islands over
these regions. Then this class is a normal family.}
\vspace{.1in}

\noindent
In fact this result was obtained by Ahlfors himself in \cite{A5} but
Dufresnoy's method permits an automatic derivation of many similar results from
Ahlfors' metric-topological theory. 

Another situation when (5.3) applies is the following. Let $D\subset\C$ be
a domain, $f:D\to D_0$ a meromorphic function
in ${\bar D}\backslash\{\infty\}$, and 
$D_0\subset\bC$ a Jordan domain. Assume that $f$ maps
$\partial D\cap\C$ to $\partial D_0$, 
Then $D$ equipped with the pullback of the
spherical metric is regularly exhaustible; one can take
an exhaustion by $D\cap\{ z:|z|\leq r\}$. It follows from Ahlfors' theory
that $f$ can omit at most two values from $D_0$ \cite[Theorem VI.9]{Tsuji}. This is how
the result about countability of the set of direct singularities for a parabolic
surface (3.3) can be proved.
The situation when $D=D_0$ (and they are not necessarily Jordan)
occurs frequently in holomorphic dynamics,
namely in the iteration theory of meromorphic functions; see, for example
\cite{Bergweiler}, where $D$ is a periodic component
of the set of normality of a meromorphic function. In this case it follows
from Ahlfors' theory that the Euler characteristic of $D$ is non-positive
\cite{Bolsch}

Here is another variation on the same topic, due to Noshiro and Kunugui,
see for example \cite{Tsuji}:
\vspace{.1in}

\noindent
{\em let $f$ be a meromorphic
function in the unit disk $\U$ such that $\lim_{z\to\zeta}f(z)$
exists and belongs to the unit circle $\partial\U$ for all
$\zeta\in\partial\U\backslash E$, where $E$ is a closed set of zero
logarithmic capacity, $|f(0)|<1$. Then $f(z)=a$ has solutions $z\in \U$ for all
$a\in\U$ with at most two exceptions. If $f$ is holomorphic
then the number of exceptional values is at most one.}
\vspace{.1in}

In \cite{A8} Ahlfors made, using his own expression, his third attempt
to penetrate the reasons behind Nevanlinna's value distribution theory.
This time he bases his investigation  on the Gauss -- Bonnet formula, which
is 
$$\int\int_RK\,d\rho=-2\pi\chi-\int_{\partial R}gds.$$
Here $R$ is a bordered surface with a smooth Riemannian metric $ds$
with associated area element $d\rho$, $K$ is the
Gaussian curvature, $\chi$ is the Euler
characteristic and $g$ is the geodesic
curvature. Choosing a metric with finitely many singularities of the type
(\ref{sing}) and pulling it back from $R_0$ to $R$ via $\pi$, one can
obtain the relation between Euler characteristics of $R_0$ and $R$.
The only analytical problem is to estimate the integrals along the boundary
curves in the Gauss--Bonnet formula. This general method permits one to give
a unified proof of
both (\ref{sch}) and the usual SMT (\ref{SMT-bis}). This paper \cite{A8}
was the basis of many generalizations of value distribution theory.
\vspace{.1in}

\noindent 
{\bf 6. Bloch constant and Ahlfors' extension of Schwarz lemma}
\setcounter{section}{6}
\setcounter{equation}{0}
\vspace{.1in}

In 1921 Valiron proved another necessary condition of
parabolic type. 
\vspace{.1in}

\noindent
{\em Let 
$$\pi:R\to\C$$
be a parabolic Riemann surface. Let $ds$ be the the pullback of the Euclidean
metric from $\C$ to $R$. Then there are disks of arbitrarily large radii
in $R$ with respect to $ds$ such that the restriction of $\pi$ to these
disks is one-to-one.} 
\vspace{.1in}

\noindent
This can be easily derived from the Five Islands theorem (5.1).
Much better known is the statement which corresponds to this via Bloch's
Principle (5.6). 
\vspace{.1in}

\noindent
{\bf Bloch's Theorem}. {\em Let $f$ be a holomorphic function in the unit
disk $\U$ satisfying $|f'(0)|=1$. Then there is a relatively compact region
$D\subset\U$ which is mapped by $f$ univalently onto a disk of radius $B$,
where $B$ is an absolute constant.}
\vspace{.1in}

\noindent
The largest value of $B$ for which this theorem is true is called Bloch's
constant. Its precise value is unknown to this day. A very plausible
candidate for the extremal
Riemann surface can be described as follows. Consider the tiling
of the plane with equilateral triangles, such that $0$ is the center of
one of the triangles. Denote the set of all vertices by $X$. There is
unique simply connected Riemann
surface $R$ which has no singular points over $\C$
and all points over $X$ are simple critical points. This $R$ is of
hyperbolic type by Valiron's theorem. Let $f$ be the
corresponding function in the unit disk (normalized as in
Bloch's theorem). In \cite{A9} 
the authors make the calculation for this function and find that
$$B\leq B':=\sqrt{\pi}\cdot 2^{1/4}
\frac{\Gamma\left(\frac{1}{3}\right)}
     {\Gamma\left(\frac{1}{4}\right)}
\left(
      \frac{\Gamma\left(\frac{11}{12}\right)}
           {\Gamma\left(\frac{1}{12}\right)}
\right)^{1/2}=0.4719\ldots.$$
This $B'$ is conjectured to be the correct value of Bloch's constant $B$.

In \cite{A10} Ahlfors gives the lower estimate $B\geq \sqrt{3}/4>0.433$, but
much more
important is the method which Ahlfors used to obtain this result.

Let $ds=p(z)|dz|$ be a conformal metric. Its Gaussian curvature is expressed 
by $-p^{-2}\Delta\log p$. This is invariant under conformal mappings.
We denote by $d\sigma$ the Poincar\'{e}
metric in the unit disk, 
$$|d\sigma|=2\frac{|dz|}{1-|z|^2},$$
whose Gaussian curvature is identically equal to $-1$. The usual conformal invariant form
of the Schwarz lemma says that every holomorphic map of the unit disk into
itself is contracting with respect to Poincar\'{e} metric.
We call $ds'=p'|dz|$ a supporting metric of $ds=p|dz|$ at the point $z_0$
if $p'(z_0)=p(z_0)$ and $p'(z)\leq p(z)$ in a neighborhood of $z_0$.
\vspace{.1in}
 
\noindent
{\bf  Theorem.} {\em Suppose that
$p$ is continuous in the unit disk and it is possible to find a supporting metric of curvature $\leq -1$ at every point.
Then} $ds\leq d\sigma$.
\vspace{.1in}

The usual Schwarz lemma follows if we take the pullback of $d\sigma$ as $ds$.
This extended form of the Schwarz lemma is frequently used now in differential
geometry and geometric function theory. See \cite{Brooks} for a survey of
these applications.

Returning to Bloch's constant, the best known lower estimate is $B>\sqrt{3}/4+
2\times 10^{-4}$ which is due to Chen and Gauthier \cite{Chen}.
\vspace{.1in}

\noindent
{\bf 7. Holomorphic curves}
\setcounter{section}{7}
\setcounter{equation}{0}
\vspace{.1in}

{}From the algebraic point of view the Riemann sphere $\bC$ is the complex
projective line ${\bf P}^1$. 
Value distribution theory can be extended to holomorphic mappings
from $\C$ to complex projective space ${\bf P}^n$. Let us recall the
definition. In $\C^{n+1}\backslash \{ 0\}$ we consider the
following equivalence relation: $(z_0,\ldots,z_n)\sim ({z'}_0,\ldots,{z'}_n)$
if there is a constant $\lambda\in \C\backslash\{ 0\}$ such that
${z'}_j=\lambda z_j$. The set of equivalence classes with the natural
complex analytic structure is 
$n$-dimensional complex projective space ${\bf P}^n$.
Let $\pi:\C^{n+1}\rightarrow{\bf P}^n$ be the projection.
Coordinates of a $\pi$-preimage of a point $w\in{\bf P}^n$
in $\C^{n+1}$ are called the homogeneous coordinates of $w$.
A holomorphic
map $f:\C\to{\bf P}^n$ is called a holomorphic curve.
It can be lifted to a map $F:\C\to\C^{n+1}$ such that $f=\pi\circ F$.
This $F=(F_0,\ldots,F_n)$ is called a homogeneous representation. Here the
$F_j$ are entire functions without common
zeros. For $n=1$ this is just a representation of a meromorphic function
as a quotient of two entire ones.
A generic point in ${\bf P}^n$ has no preimages under $f$ so one studies
preimages of hypersurfaces, in particular, preimages of hyperplanes.
An analog of Nevanlinna's theory for this case was created by H. Cartan
\cite{Cartan}, and the role of the spherical metric is now played by a Hermitian
positive $(1,1)$-form, the so-called Fubini-Study form. The pullback of this 
form via a holomorphic curve $f$ is
$$d\rho=\| F\|^{-4}\sum_{i>j}|F_i{F'}_j-{F'}_iF_j|^2\,dm,$$
where $\|.\|$ is the usual Euclidean norm in $\C^{n+1}$ and $dm$ is the
area element in $\C$. The Nevanlinna-Cartan characteristic of $f$ is defined
by
$$T_f(r)=\int_0^t\frac{A_f(t)}{t}dt,\quad\mbox{where}\quad A(r)=\int_{|z|\leq r}d\rho.$$
A hyperplane $A\subset{\bf P}^n$ is given by a linear equation in homogeneous
coordinates: $a_0z_0+\ldots+a_nz_n=0$, so the preimages of this hyperplane
under $f$ are just zeros of the linear combination $g_A=a_0F_0+\ldots+a_nF_n$.
The number of these zeros in the disk $|z|\leq r$ is denoted by $n(r,A)$ and
the averaged counting functions $N(r,A)$ are defined as before. The proximity
functions $m(r,A)$ can be also defined by analogy with the one-dimensional
case, the role of the chordal distance now being played by the sine of
the angle between a line and a hyperplane:
$$m(r,A)=\varint\log\frac{\| F\|\cdot\| A\|}{| g_A|}\,d\theta.$$
The First Main Theorem of Cartan says that
\begin{equation}
\label{C-FMT}
m(r,A)+N(r,A)=T(r)+O(1)\quad\mbox{\em for every hyperplane}\; A,
\end{equation}
where we assume only that $f(\C)$ is not contained in $A$.
In the Second Main Theorem of Cartan we assume that several hyperplanes
$A_1,\ldots,A_q$ are in {\em general position} that is every $n+1$ of them
have empty intersection. We also denote by $N^*(r)$ the averaged counting
function of zeros of the Wronskian determinant $W(F_0,\ldots,F_n)$; these 
zeros do not depend on the choice of homogeneous representation $F$. To guarantee
that this Wronskian is not identical zero we need to assume that our curve
is {\em linearly non-degenerate} that is $f(\C)$ is not contained in any
hyperplane, which is the same as to say that $F_0,\ldots,F_n$ are linearly
independent.
The SMT of Cartan is
\begin{equation}
\label{C-SMT}
\sum_{j=1}^qm(r,A_j)+N^*(r)\leq (n+1) T(r)+S(r),
\end{equation}
where the error term $S(r)$ has the same estimate as in one-dimensional
Nevanlinna theory. 

One important corollary is the Borel theorem: {\em a curve omitting
$n+2$ hyperplanes in general position has to be linearly degenerate}.
This is not a ``genuine'' Picard-type theorem, because the conclusion is
``linearly degenerate'' rather than ``constant''. But one can deduce from
Borel's theorem the following: {\em a curve omitting $2n+1$ hyperplanes in general
position has to be constant}. This corollary from Borel's theorem
was apparently formulated for the first
time by P. Montel \cite{Montel}. 

It is strange that Cartan's paper \cite{Cartan} was overlooked by many
researchers 
for about 40 years. 
Even as late as in 1975 some
specialists in holomorphic curves were surprised to hear about Cartan's work.

Meanwhile H. Weyl and his son J. Weyl in 1938 restarted the subject of value
distribution of holomorphic curves from the very beginning \cite{Weyls},
independently of Cartan.
They failed
to come to conclusive results like Cartan's second main theorem, but they
contributed two important ideas to the subject. The first of them
was to consider the so-called {\em associated curves}.
Here we need some definitions. A set of $k+1$ linearly independent vectors
$A_0,\ldots,A_k\in\C^{n+1}$ defines a $k+1$-dimensional
linear subspace. It is convenient to use the wedge product $A_0\wedge\ldots
\wedge A_k$ to describe this subspace. Two $(k+1)$-tuples define the
same subspace if and only if their wedge products are proportional. 
The elements of $\wedge^{k+1}\C^{n+1}$ are called polyvectors; those of them
which have the form $A_0\wedge\ldots\wedge A_k$, where $A_j\in \C^{n+1}$
are called decomposable polyvectors.
If we fix a basis in $\C^{n+1}$ this defines a basis in $\wedge^{k+1}\C^{n+1}$.
Thus to each $k+1$-subspace corresponds a point in  
${\bf P}^{N_k}$, where $N_k=\binom{n+1}{k+1}-1$, the homogeneous coordinates
of this point being coordinates of $A_0\wedge\ldots\wedge A_k$ in the fixed basis
of $\wedge^{k+1}\C^{n+1}$. The coordinates we just introduced are called the
Pl\"{u}cker coordinates. 
A linear subspace of dimension $k+1$ in $\C^{n+1}$ projects into a
projective subspace of dimension $k$. The set of all such projective subspaces
identify $G(k,n)$ with a manifold in 
${\bf P}^{N_k}$ of dimension $(k+1)(n-k)$.
An explicit system of equations defining this manifold
can be written (see for example \cite{HP,GK}).

To every holomorphic curve $f$ in ${\bf P}^n$, the Weyls  associate curves
$f_k:\C\to G(n,k)$ which in homogeneous coordinates have the representation
$F\wedge F^{(1)}\wedge F^{(2)}\wedge\ldots\wedge F^{(k)},\;
 k=0,\ldots, n-1;$ where $F^{(j)}$
is the $j$-th derivative. The geometric
interpretation is that they assign to each point $f(z)$ on the curve
the tangent line $f_1(z)$ and the {\em osculating $k$-planes} $f_k(z)$.
For each $k$,  $f_k(z)$ is the unique $k$-dimensional projective
subspace which has a
contact of order at least $k$ with the curve at the point $f(z)$.

The second important idea of the Weyls was to consider the proximity functions
of a curve $f$ with respect to projective subspaces of any codimension, not only
with respect to hyperplanes. Such a proximity function may be obtained by
averaging the hyperplane proximity functions over all hyperplanes
containing the given subspace. The FMT (\ref{C-FMT}) no longer holds for
$k$-subspaces with $k<n-1$ because they normally have no preimages at all.
But the SMT interpreted as an upper bound for proximity functions still
has a sense. The Weyls \cite{Weyls} managed to prove such SMT for points, which
are
$0$-subspaces. 

It was Ahlfors who proved in \cite{A11} the precise estimates for
proximity functions for all $f_k$ in all codimensions.
To formulate the result we use the following notation.
Let $A^h$ be a decomposable $h$-polyvector. 
Then
$$m_k(r, A^h)=\varint\log\frac{\| F_k\|\cdot\| A^h\|}{\| F_k\vee A^h\|}d\theta.$$
The expression in the logarithm is reciprocal to
the distance between the
polyvector $F_k$ and the subspace orthogonal to $A_h$. It is expressed in
terms of the inner product $\vee$ of polyvectors, which coincides with the usual
dot product when $h=k$. 
The number of critical points of the map $f_k:\C\to
G(n,k)$ is denoted by $n_k^*(r)$ and the corresponding averaged counting
function by $N_k^*(r)$.  We also set $T_k=T_{f_k}$.
Now we assume that $f$ is a linearly non-degenerate
curve and a finite system of decomposable
polyvectors in general position is given in each dimension $h$.
Ahlfors' result is the following
(we use the formulation from \cite{Wu}).
\newline
{\em For $0\leq h\leq k$
$$\sum_{A^h}m_k(A^h)\leq\binom{n+1}{h+1}T_k-
\sum_{m=k}^{n-1}\sum_{i=0}^h p_h(m-i,h-i)N_{m-i}^*-$$
$$-\sum_{i=k-h-1}^k
\binom{i}{h-k+i+1}\binom{n-i-1}{k-i}T_i+S(r),
$$
and for $k\leq h\leq n-1$
$$\sum_{A^h}m_k(A^h)\leq\binom{n+1}{h+1}T_k-
\sum_{m=0}^k\sum_{i=0}^{n-k-1}p_{n-h-1}(n-m-i-1,m-h-i-1)N_{m+i}^*-$$
$$-\sum_{i=k}^{n-h+k}
\binom{n-1-i}{h-k-1}\binom{i}{k}T_i+S(r),
$$
where the following notation is used
$$p_h(k,l)=\binom{n+1}{h+1}-
\sum_{j\geq 0}\binom{k+1}{k+j+1}\binom{n-k}{h-l-j}\geq 0
.$$
Here the binomial coefficient is defined for all integers by} $(1+x)^n=\sum_k
\binom{n}{k}x^k$.

In particular for $h=k$, when the First Main Theorem applies
we have
\begin{equation}
\label{hyperplanes}
\sum_{A^k}m_k(A^k)\leq\binom{n+1}{k+1}T_k-
\sum_{m=0}^k\sum_{i=m}^{n+m-k-1}p_{k}(i,m)N_i^*+S(r).
\end{equation}
{}From this one can derive
\begin{equation}
\label{hyperplanes2}
\sum_{A^k}m_k(A^k)\leq\binom{n+1}{k+1}T_k
\end{equation}
This has the same form as Cartan's SMT (\ref{C-SMT}) would give if one
applies it to the associated curves $f_k$ and drops the $N^*$ term. 
But (\ref{hyperplanes2}) {\em does not} follow from Cartan's
SMT. The catch is that even if $f$ is linearly non-degenerate, the associated
curves $f_k$ might be linearly degenerate. One can also show (see, for example
\cite[p. 138]{Fujimoto}) that the ramification term $N^*$ in Cartan's SMT
can be obtained from (\ref{hyperplanes}). Thus Ahlfors' result is stronger
then Cartan's SMT when applied to the associated curves $f_k$, and decomposable
hyperplanes, and they coincide
when applied to the curve $f$ itself and hyperplanes.

The difficulties Ahlfors had to overcome in this paper are enormous.
The main idea can be traced back to his proofs of Nevanlinna's SMT,
but multidimensionality causes really hard problems. As Cowen and Griffiths 
say in \cite{CG} ``The Ahlfors theorem strikes us as one of the few instances
where {\em higher codimension} has been dealt with {\em globally} in
complex-analytic geometry''.
Ahlfors in his commentaries \cite[p. 363]{CO} says: ``In my own eyes the paper
was one of my best, and I was disappointed that years went by without signs
that it had caught on". 

Since he wrote this the situation has changed a little. H. Wu \cite{Wu}
published a detailed and self-contained exposition of Ahlfors' work.
Two new simplified  versions of the proof of (\ref{hyperplanes}) were
given in 
\cite{CG} and \cite{Fujimoto}; 
the second work gives important applications to
minimal surfaces in ${\bf R}^n$. Still, from my point of view much
in this work of Ahlfors
remains unexplored. For example, let $f$ be a holomorphic curve
in ${\bf P}^2$ and $a_1,\ldots,a_q$ a collection of points in ${\bf P}^2$.
What is the smallest constant $K$ such that
\begin{equation}
\label{conjecture}
\sum_{j=1}^qm_f(r,a_j)\leq KT_f(r)+S(r)
\end{equation}
holds? Ahlfors' relations will give this with $K=3/2$; the same can be 
deduced from (\ref{C-SMT}). But one can conjecture that in fact $K=1$.
A related question was asked by Shiffman in \cite{Shiffman}; no progress
has been made. 
\vspace{.1in}

\noindent
{\bf 8. Multi-dimensional counterpart of the type problem.}
\setcounter{section}{8}
\setcounter{equation}{0}
\vspace{.1in}

Here we will mention very briefly some later research in the spirit of
the type problem where Ahlfors' ideas play an important role.
\vspace{.1in}

{\em Quasiregular mappings of Riemannian manifolds}. A mapping $f$ between two
$n$-dimensional Riemannian manifolds is called $K$-quasiregular
if it belongs to the Sobolev class $W^1_{n,{\rm loc}}$ and its derivative
almost everywhere satisfies $|df|^n\leq KJ_f$, where $J_f$ is the Jacobian.
Given two orientable Riemannian manifolds $V_1$ and $V_2$, one may ask whether
a non-constant quasiregular map $f:V_1\to V_2$ exists. To answer this question
Gromov \cite[Ch. 6]{Gromov} generalizes the parabolic type criterion from Section 1
in the following way. Let $d(V)$ be the supremum of real numbers $m>0$ such that
for some constant $C$ the isoperimetric inequality $\mbox{volume}(D)
\leq C\mbox{area}(\partial D)^{m/(m-1)}$ holds for every compact $D\subset V$.
Here ``area'' stands for the $n-1$ dimensional measure. This number $d(V)$
is called the {\em isoperimetric dimension} of $V$. The isoperimetric dimension of 
${\bf R}^n$ with the Euclidean metric is $n$ and the isoperimetric dimension
of the hyperbolic (Lobachevskii) space ${\bf H}^n$ is $\infty$. 
Ahlfors' argument presented in Section 1 proves the following:
\vspace{.1in}

\noindent
{\em Assume that for some point $a\in V$ we have
\begin{equation}
\label{ahlfors}
\int^\infty\frac{dr}{L^{1/(m-1)}(r)}=\infty,
\end{equation}
where $L(r)$ is the $(n-1)$-measure of the sphere of radius $r$ centered at $a$.  
Then the isoperimetric dimension of any conformal metric on $V$ is at most $m$}.
\vspace{.1in}

As a corollary, one deduces that a nonconstant quasiregular map from $V_1$
to $V_2$, where $V_1$ satisfies (\ref{ahlfors}), is possible only if the isoperimetric
dimension of $V_2$ is at most $m$. 

Because $f$ can be lifted to a map ${\tilde f}:V_1\to{\tilde V_2}$ to the
universal covering, one is interested in estimating $d({\tilde V}_2)$.
Gromov proves that for compact $V_2$ the isoperimetric dimension of ${\tilde V_2}$ depends only on the fundamental group: a distance and volume can be
introduced on every finitely generated group, and the fundamental group will
satisfy the same isoperimetric inequality as~${\tilde V}_2$.

As a simple corollary one obtains
\vspace{.1in}

\noindent\addtocounter{equation}{+1}
(\theequation)$\;$ {\em There is no quasiregular
map from ${\bf R}^3$ to
$S^3\backslash N$, where $S^3$ is the $3$-dimensional sphere 
and
$N$ is a non-trivial knot.}
\vspace{.1in}

This result is superseded by the following generalization of
Picard's theorem proved by Rickman in 1980
\vspace{.1in}

\noindent\addtocounter{equation}{+1}
(\theequation) {\em For every $n\geq 2$ and $K\geq 1$ there exists $q$
 such that a $K$-quasiregular map ${\bf R}^n\to S^3$ cannot omit more
then $q$ points.}
\vspace{.1in}

Holopainen and Rickman \cite{HR} later proved that $S^3$ can be replaced with 
any compact manifold. The simplest proof of (8.3) is due to Lewis \cite{Lewis}.
It is also one of the most elementary proofs of Picard's
theorem.

The experts were stunned when Rickman constructed examples \cite{Rick}
showing that the number of omitted values may really depend on $K$. Any
finite number of values can be actually omitted.
\vspace{.1in}


\end{document}